\documentclass[11pt, a4paper]{article}  
\author{Frank Schuhmacher}

\title{Noncommutative complex analytic spaces}

\usepackage{mathenv}

\usepackage[intlimits]{amsmath}

\usepackage[T1]{fontenc}

\usepackage{amssymb}

\usepackage{amsthm}

\usepackage{amscd}

\usepackage[all]{xypic}

\usepackage{makeidx}

\newcounter{punkt}

\theoremstyle{definition}

\newtheorem{defi}{Definition}[section]

\newtheorem{beisp}[defi]{Example}

\theoremstyle{theorem}

\newtheorem{satz}[defi]{Theorem}

\newtheorem{kor}[defi]{Corollary}

\newtheorem{prop}[defi]{Proposition}

\newcommand{\nach}{\longrightarrow}

\newcommand{\sub}{\subseteq}

\newcommand{\isom}{\cong}

\newcommand{\RR}{ \mathbb{R} }

\newcommand{\ZZ}{\mathbb{Z}}

\newcommand{\CC}{\mathbb{C}}

\newcommand{\Oh}{\mathcal{O}}

\newcommand{\ppp}{\cdot\ldots\cdot}

\newcommand{\K}{\mathcal{K}}

\newcommand{\I}{\mathcal{I}}

\newcommand{\J}{\operatorname{J}}

\newcommand{\m}{\mathfrak{m}}

\newcommand{\red}{\operatorname{red}}

\newcommand{\ddd}{,\ldots,}

\newcommand{\kkk}{+\ldots+}

\newcommand{\bbb}{|\ldots|}
\newcommand{\ab}{\operatorname{ab}}
\newcommand{\abs}{\operatorname{abs}}
\newcommand{\unab}{\operatorname{unab}}
\newcommand{\op}{\operatorname{op}}
\newcommand{\ec}{\hat{e}}

\newcommand{\alphac}{\hat{\alpha}}
\newcommand{\otc}{\hat{\otimes}}

\begin{document}
\maketitle
\begin{abstract}
In this paper, we define
NC complex spaces as 
complex spaces together with a structure sheaf of  
associative algebras in such a way that the abelization of the
structure sheaf is the sheaf of holomorphic functions.
\end{abstract}

\section*{Introduction}
We want to
study noncommutative structures
on underlying classical geometric objects.
A few words to motivate this objective:
Let (spaces) be one of the following    
categories: (schemes), (varieties),
(analytic spaces), (smooth manifolds) or an appropriate subcategory
of one of them. In particular, (spaces) is a subcategory of
the category of locally ringed spaces, i.e. a space
consists of a topological space
$X$ and a structure sheaf $\Oh$ of commutative rings and perhaps
of extra data.
Deformation quantization leads to the
problem of finding an extension (NC spaces) of the category
(spaces), such that for a space $(X,\Oh)$ and a non commutatively
deformed structure sheaf $\Oh'$ with abelization $\Oh$, the pair
$(X,\Oh')$ still belongs to (NC spaces). More generally, if
$\Oh'$ is a fibre of a deformation of the structure sheaf
$\Oh$ of $X$, there should
should be an object $(X',\Oh')$ of (NC spaces), where $X'$ is obtained by a
classical deformation of $X$.\\

Kapranov \cite{Kapr} has defined an extension 
(NC schemes) of the category (schemes) in which at least formal 
noncommutative deformations of the structure sheaf of a scheme live.
In this paper, we construct such an extension
(NC complex analytic spaces) of the category of complex analytic 
spaces. An advantage of the analytic context is that
we can allow \textit{local}
noncommutative deformations, i.e. deformations
parametrized by a space germ, not only \textit{formal} deformations,
within the category (NC complex analytic spaces). 
This paper only contains the definition and some examples
of NC complex analytic spaces. Their deformations will be described
in a forthcomming paper.\\

In the first section, we give the prerequisites on noncommutative
power series and introduce ``convergent'' noncommutative power
series. In the second section, we construct
a natural noncommutative structure sheaf on the n-dimensional
affine complex space and show that local models of analytic
spaces can be equipped with a noncommutative structure sheaf.  
NC complex analytic spaces are just defined as the objects obtained
by glueing such NC local models.
In the third section, we give several examples of NC complex
analytic spaces, like supermanifolds and NC projective 
varieties.

\section{Noncommutative power series}

The starting point of the theory of NC complex spaces
is to replace power series where-ever they occur in the
theory of holomorphic functions by noncommutative power
series. 

\paragraph{Multi-indices}
For each multi-index $I=(i_1\ddd i_k)$ in $\{1\ddd n\}^k$, set
\begin{align*}
x^I&:=x_{i_1}\ppp x_{i_k},\\
\#I&:=k.
\end{align*}
By definition, there is an empty multi-index $0$ 
with $x^0=1$, $\#0=0$. 
For multi-indice $I,J$, set $I+J:=(i_1\ddd i_{\#I},j_1\ddd j_{\#J})$.  
Let $I<J$ be the relation
$\#I<\#J$ and $I$ is obtained from $J$ by deletion of indices.
Observe that, if $I\leq J$ (i.e. $I<J$ or $I=J$), then there exists
an order-preserving, injective map $\alpha$ from
the set $\{1\ddd\#I\}$ into $\{1\ddd\#J\}$ such that
$i_k=j_{\alpha(k)}$, for $k=1\ddd\#I$. 
We say that $\alpha$ is a map from $I$ to $J$.
Write
${J \choose I}$ for the number of such maps $\alpha:I\nach J$. 
By definition, ${J \choose 0}=1$, for any $J$.
Obviously, ${J \choose I}\leq {\#J\choose\#I}$.
Set $J-_\alpha I$ to be the multi-index obtained from $J$ by deletion
of $j_{\alpha(1)}\ddd j_{\alpha(\#I)}$. For elements $p\in\CC^n$,
we can just write $p^{J-I}$ insetad of $p^{J-_\alpha I}$, since
this value does not depend on the map $\alpha$.
A multi-index $I$ is called \textbf{ordered}, if
$i_k\leq i_l$, for $k<l$.

\paragraph{Formal and convergent power series}
For an element $p=(p_1\ddd p_n)\in\CC^n$, we 
write $\CC[[x_1-p_1\bbb x_n-p_n]]$ for the noncommutative
power series algebra in formal variables $x_i-p_i$ and
$\CC[x_1-p_1\bbb x_n-p_n]$ for the subalgebra of noncommutative
polynomials. We have a canonical epimorphism $\ab$
from $\CC[[x_1-p_1\bbb x_n-p_n]]$
to the commutative power series algebra $\CC[[x_1-p_1\ddd x_n-p_n]]$.
Set $f_{\ab}:=\ab(f)$.
A noncommutative power series 
$$f=\sum_Ia_I(x-p)^I=\sum_{m=0}^\infty\sum_{\#I=m}a_I(x-p)^I$$
in $\CC[[x_1-p_1\bbb x_n-p_n]]$
is called \textbf{convergent}, 
if the image $f_{\abs}$ of $\sum_I|a_I|(x-p)^I$
under $\ab$ is convergent. Write $\CC\{x_1-p_1\bbb x_n-p_n\}$
for the algebra of convergent noncommutative power series, and
for $r\in\RR_{>0}^n$, write
$\CC_r\{x_1-p_1\bbb x_n-p_n\}$ for the subalgebra consisting of those
$f$ such that $f_{\abs}$ converges on the open polydisk
$P(p,r)$ with multiradius $r$ centered at $p$. Given such an $f$, for
elements $q$ in $P(p,r)$, set $f(q):=f_{ab}(q)$.
Lets make the convention that $\CC\{x-p\}$ stands
for $\CC\{x_1-p_1\bbb x_n-p_n\}$ (always with $n$ variables)
and that $\CC\{y-q\}$ stands
for $\CC\{y_1-q_1\bbb y_m-q_m\}$ (always with $m$ variables).

A noncommutative power series
$f=\sum_Ia_I(x-p)^I$
is called \textbf{commlike}, if $a_I=0$, for nonordered $I$.
The map $\ab$ from the noncommutative to
the commutative power series ring has a splitting $\unab$
in
the category of $\CC$-vectorspaces, such that the composition
$\unab\circ\ab$ fixes exactly the commlike
power series.

\paragraph{Morphisms}
If $g=\sum_Jb_J(x-p)^J$ is a power series in $\CC[[x-p]]$
and if the power series
$f_1\ddd f_m$ with $f_\nu=\sum_Ia_{\nu,I}(y-q)^I$
belong to the maximal ideal $\m_{[[y-q]]}$
of
$\CC[[y-q]]$, we can form the power series
$$g(f_1,...,f_m):=\sum_J(\sum b_K\cdot
a_{k_1,I_1}\ppp a_{k_{k_{\#K},I_{k,\#K}}})(y-q)^J,$$
where the sum in the bracket is over all multi-indices
$K$ such that $J=I_{k_1}\kkk I_{k_{\#K}}$.

A \textbf{morphism} $\CC[[x-p]]\nach\CC[[y-q]]$
of noncommutative power series algebras is a local algebra
homomorphism of the form $g\mapsto g(f_1\ddd f_n)$,
for a given $n$-tuple $(f_1\ddd f_n)$ of elements of
$\m_{[[y-q]]}$.

\begin{prop}\label{composition}
If $g$ and each $f_i$ are convergent power series, then
$g(f_1\ddd f_m)$ is also convergent. 
\end{prop}

A \textbf{morphism} $\CC\{x-p\}\nach\CC\{y-q\} $
of convergent noncommutative power series algebras
is a local algebra homomorphism of the form
$g\mapsto g(f_1\ddd f_n)$, for a given $n$-tuple
$(f_1\ddd f_n)$ of elements of $\m_{\{y-q\}}$. 
\begin{prop}
Let $f:\CC[[x-p]]\nach\CC[[y-q]]$ and $g:\CC[[y-q]]\nach\CC[[z-r]]$
be given by $f(x_\nu-p_\nu)=\sum_Ia_{\nu,I}(y-q)^I$ and
$g(y_\mu-q_\mu)=\sum_{J}b_{\mu,J}(z-r)^J$. Then, we have
$g(f(x_\nu-p_\nu))=\sum_Kc_K(z-r)^K$, with
$$c_K=\sum_Ia_{\nu,I}\sum b_{i_1,J_1}\ppp b_{i_{\#I},J_{\#I}},$$
where the second sum is taken over all multi-indices
$J_1\ddd J_{\#I}$ such that $J_1\kkk J_{\#I}=K$.
\end{prop}

For an endomorphism $f$ of
$\CC[[x-p]]$, with $f(x_\nu-p_\nu)=\sum_Ia_{\nu,I}(x-p)^I$, set
$Jf$ to be the $n\times n$-matrix $(a_{\nu,i})_{\nu,i}$.
\begin{satz}
An endomorphism $f$ of $\CC[[x-p]]$ is an automorphism, if and
only if $Jf$ is invertible.
\end{satz}
\begin{proof}
Suppose that $Jf$ is invertible. Inductively, for
$k\geq 1$, $\mu=1\ddd n$ and $J\in\{1\ddd n\}^\mu$,
we  will define coefficients $b_{\mu,J}\in\CC$
such that the endomorphism $g^{(k)}$ of $\CC[[x-p]]$ with
$g^{(k)}(x_\mu-p_\mu)=\sum_{\#J\leq k}b_{\mu,J}(x-p)^J$ is inverse to $f$
modulo $\m_{[[x-p]]}^{k+1}$.
For $k=1$, by Proposition~\ref{composition}, the necessary
(and sufficient) conditions are
\begin{align*}
\sum_{i=1}^na_{\nu,i}b_{i,\nu}=&1\quad\text{ for all } \nu,\\
\sum_{1=1}^na_{\nu,i}b_{i,\mu}=&0\quad\text{ for all } \nu\neq\mu.
\end{align*}
We can find such coefficients, if and only if $Jf$ is invertible.
Now suppose that $Jf$ is invertible and that the coefficients
$b_{\mu,J}$ are constructed adequately, for $\#J\leq k-1$.
For fixed $J$ with $\#J=k$, we have to
find $b_{\mu,J}$ such that 
$$\sum_{i=1}^na_{\nu,i}b_{i,J}+\sum_{\#I\geq 2}a_{\mu,I}
\sum b_{i_1,J_1}\ppp b_{i_{\#I},J_{\#I}}=0.$$
The second sum is known, and since $Jf$ is invertible, we
can find adequate $b_{1,J}\ddd b_{n,J}$.
\end{proof}

\begin{kor}
Each lift of an automorphism of a commuative (convergent)
power series ring to an endomorphism of the noncommutative
(convergent) power series ring is again an automorphism. 
\end{kor}

\paragraph{Complete tensor products}
For power series algebras $\CC[[x]]=\CC[[x_1|...|x_n]]$ and 
$\CC[[y]]=\CC[[y_1|...|y_m]]$, we
define the \textbf{free product} 
$$\CC[[x]]\ast \CC[[y]]:=\CC[[x_1|...|x_n|y_1|...|y_n]]$$
and the \textbf{analytic tensor product} 
$\CC[[x]]\otc_\CC\CC[[y]]$ as the power series
algebra in $x_1,...,x_n$ and $y_1,...,y_m$, where the $y_j$ are 
assumed to commute with the $x_i$.
We make the corresponding definitions for convergent power series
algebras. 

\paragraph{Finitely generated ideals}
In the NC analytic context, the concept of finitely generated
two-sided ideals has to be slightly adapted. As a reason,
we give the following example:
\begin{beisp}
Let $(x)$ be the two-sided ideal of $\CC[[x|y]]$,
consisting of all noncommutative power series
where at least one factor $x$ arises in each monomial.
Observe that $(x)$ is not the two-sided ideal generated by $x$ in
the algebraic sense.
\end{beisp} 
\begin{proof}
Assume that $(x)$ is the two-sided ideal generated by $x$ in the 
algebraic sense. Then we can find power series 
$f_i=\sum_j a_{ij}y^j$ and $g_i=\sum_jb_{ij}y^j$ in
$\CC[[y]]$, $i=1\ddd N$ such that
$$xyx+y^2xy^2+...=\sum_{i=1}^Nf_i(y)\cdot x\cdot g_i(y).$$
The right hand-side takes the form
$\sum_{j,k}(\sum_i a_{ij}b_{ik})y^jxy^k$.
In particular, for $j,k\leq N+1$, we would get
$\sum_i^Na_{ij}b_{ik}=\delta_{j,k}$,
which is impossible, since the left hand-side is
a product of two matrices of rank at most $N$.
\end{proof}

By definition, the opposite algebra 
$\CC[[x]]^{\op}$ $\{f^{\op}:\;f\in\CC[[x]]\}$ of $\CC[[x]]$
is the set $\{f^{\op}:\;f\in\CC[[x]]\}$ with operations
$f^{\op}_1+f^{\op}_2=(f_1+f_2)^{\op}$ and 
$f_1^{\op}\cdot f_2^{\op}=(f_2\cdot f_1)^{\op}$. Observe that
$\CC[[x]]^{\op}$ is naturally isomorphic to the power
series algebra $\CC[[x_1^{\op}|...|x_n^{\op}]]$ and the assignment
\newcommand{\OP}{\operatorname{OP}}
$x\mapsto x^{\op}$ defines an isomorphism 
$\OP:\CC[[x]]\nach \CC[[x]]^{\op}$.
Attention, in general, $\OP(f)\neq f^{\op}$, for example, 
$\OP(x_1 x_2)=x_1^{\op} x_2^{\op}=(x_2 x_1)^{\op}$.
We define the \textbf{(complete) envelopping algebra} 
$\CC[[x]]^{\ec}$ of $\CC[[x]]$
as $\CC[[x]]\otc_\CC\CC[[x]]^{\op}$.
Consider the natural epimorphism
$$\alphac:\CC[[x]]^{\ec}\nach \CC[[x]].$$
For a two-sided ideal 
$J\sub \CC[[x]]$, the inverse image $\alphac^{-1}(J)$ is not, in general
a left ideal of $\CC[[x]]^{\ec}$, since it is not, in general, 
closed under left multiplication by elements of $\CC[[x]]^{\ec}$.
We define the \textbf{completion} $\hat{J}$ of $J$ as the image under
$\alphac$ of the $\CC[[x]]^{\ec}$-left ideal generated by 
$\alphac^{-1}(J)$.
For simplicity, for elements $f_1,...,f_m$ in $\CC[[x]]$, we shall write
$(f_1,...,f_m)$ for the completion of the
two-sided ideal generated by the $f_i$. A two-sided ideal of
the form $(f_1,...,f_m)$ will be called \textbf{finitely generated}.

\begin{prop}
The Kernel $\K$ of the abelization 
$\ab:\CC[[x_1|...|x_n]]\nach\CC[[x_1\ddd x_n]]$ is
finitely generated by the commutators $[x_i,x_j]=x_ix_j-x_jx_i$,
for $1\leq i<j\leq n$.
\end{prop}
\begin{proof}
We show that for each noncommutative power series $f$,
the difference $f-f_{\ab,\unab}$ is in the image
under $\alphac$ of the $\CC[x]^{\ec}$-left ideal
generated by the commutators $[x_i,x_j]$, for $i<j$.
Without restriction, let
$f$ be the sum of its homogeneous components
$f_k$ of degree $k\geq 2$.
Each difference $f_k-f_{\ab,\unab,k}$ is of the
form $\alphac(\sum_{i<j}c^{ij}_k\cdot [x_i,x_j])$, for
certain homogeneous $c_k^{i,j}$ in $\CC[x]^{\ec}$ of
degree $k-2$.
Thus $f-f_{\ab,\unab}$ is the image under $\alphac$ of
$\sum_{i<j}(\sum_k c_k^{ij})\cdot[x_i,x_j]$.  
\end{proof}

\paragraph{Locality}
A family $(h_\alpha)_{\alpha\in A}$ of power series in $\CC[[x-p]]$
is called \textbf{summable}, if, for each multi-index
$I$, there are only finitely many $\alpha\in A$ such that
$h_{\alpha,I}\neq 0$. In this case, we can form the power
series
$$\sum_{\alpha\in A}h_\alpha:=\sum_I(\sum_\alpha h_{\alpha,I})(x-p)^I.$$

\begin{prop}\label{localring}
Both, $\CC[[x-p]]$ and $\CC\{x-p\}$, are local rings with maximal
two-sided ideal $\m_p$ generated by $x_1-p_1\ddd x_n-p_n$.
\end{prop}
\begin{proof}
It suffices to show that each element $f$ of $\CC[[x-p]]\setminus\m_p$
is a left-and right unit. Without restriction, say $f_0=1$.
Then the family $(1-f)^j;j\geq 0$ is summable.
We have 
\begin{align*}
f\cdot\sum_{j=0}^\infty(1-f)^j&=(1-(1-f))\sum_{j=0}^\infty(1-f)^j=\\
&=\sum_{j=0}^\infty(1-f)^j-\sum_{j=1}^\infty(1-f)^j=1.
\end{align*}
Thus $f$ is a left unit. 
In the same way, we show that $f$ is a right unit.
If $f$ is convergent, the sum $\sum(1-f)^j$ is also convergent. This
follows exactly as in the commutative case.
\end{proof}

\section{NC complex analytic spaces}

\paragraph{NC functions on polydisks}

\begin{prop}\label{dieerste}
Fix a point $p$ in $\CC^n$.
If $f=\sum_Ia_I(x-p)^I$ is a noncommutative
power series in $\CC_r\{x-p\}$ converging on the open
polydisk $P(p,r)$ of multi-radius $r$,
then, for each muti-index $J$ with $J\leq I$, the 
power series
$$a_J(z):=\sum_{\alpha} a_I (z-p)^{I-_\alpha J},$$
where the sum is taken over all maps $\alpha:J\nach I$ of multi-indices,
converges on $P(p,r)$.
\end{prop}
\begin{proof}
We must show that, for $q\in P(p,r)$, the series
$$a_J(q)=\sum_I{I\choose J} a_I (q-p)^{I-J}$$
converges absolutely. We may assume that
all $a_I$ are positive real numbers and that
$(p-q)$ is in $\RR^n_{\geq 0}$. The product $a_J(q)\cdot(q-p)^J$ 
lower or equals
$\sum_{I\geq J}{\#I\choose\#J}a_I(q-p)^I$, which lower
or equals
$\sum_m^\infty m^{\#J}\sum_{\#I=m}a_I(q-p)^I$.
The latter converges by the quotient criterium, since
we know that $f$ converges absolutely at $q$.
\end{proof}

The following statement is an immediate consequence
of the correponding statement for commutative power
series: 

\begin{prop}\label{diezweite}
Let $r'\in\RR_{>0}^n$ be a mult-radius such that the
open polydisk $P(q,r')$ is contained in $P(p,r)$.
The power series $f_q:=\sum_Ja_J(q)(z-q)^J$ converges
on $P(q,r')$ and represents the function $f_{\ab}|_{P(q,r')}$.
\end{prop}

\begin{prop}\label{alpha}
Let $p$ and $q$ be points in $\CC^n$ and 
$r,r'\in\RR_{>0}^n$ multi-radius
such that $P(q,r')$ is contained in the open polydisk $P(p,r)$.
The map
\begin{align*}
\alpha_r(p,q):\CC_r\{x-p\}&\nach\CC_{r'}\{x-q\}\\ 
\sum_Ia_I(x-p)^I&\mapsto\sum_Ja_J(q)(x-q)^J,
\end{align*}
with $a_J(q)$ as in Proposition~\ref{dieerste},
is an algebra homomorphism.
\end{prop} 
\begin{proof}
Consider two elements $f=\sum_Ia_I(x-p)^I$
and $\tilde f=\sum_I\tilde a_I(x-p)^I$ of $\CC_r\{x-p\}$.
We have 
$$f\cdot\tilde f=\sum_K(\sum a_I\cdot \tilde a_{\tilde I})(x-p)^K,$$
where the sum $b_K$ in the bracket is over all $I$ and $\tilde I$ such
that $K=(I,\tilde I)$.
We get
$$\alpha_r(p,q)(f\cdot\tilde f)=\sum_{L}b_L(q)(x-q)^L,$$
where $b_L(q)$ is the sum
$$\sum_{K\geq L}{K \choose L}b_K(q-p)^{K-L}=
\sum_{K\geq L}\sum{K \choose L}a_I\cdot \tilde a_{\tilde I}(q-p)^{K-L}.$$
For a given $L$, the 4-tuples $(K,\alpha,I,\tilde I)$, where
$\alpha:L\nach K$ is a map of multi-indices and $I,\tilde I$ are
multi-indices such that $K=(I,\tilde I)$, are in one-to-one correspondence
with the 6-tuples $(J,\tilde J,I,\tilde I,\alpha_1,\alpha_2)$,
where $J,\tilde J$ are multi-indices such that $L=(J,\tilde J)$ and
$\alpha_1:\J\nach I$ and $\alpha_2:\tilde J\nach\tilde I$
are maps of multi-indices.
To the 4-tuple $(K,\alpha,I,\tilde I)$, assign 
$(\alpha^{-1}(I),\alpha^{-1}(\tilde I),I,\tilde I,\alpha|_J,
\alpha|_{\tilde J})$. To the 6-tuple $(J,\tilde J,I,\tilde I,\alpha_1,
\alpha_2)$, assign $((I,\tilde I),(\alpha_1,\alpha_2),I,\tilde I)$.
Observe that, for corresponding tuples, we have
${K\choose L}={I\choose J}\cdot{\tilde I\choose\tilde J}$.
We can thus write $\alpha_r(p,q)(f\cdot\tilde f)$ as
$\sum_Lc_L(x-q)^L$, where
$$c_L=\sum\sum_{I\geq J}\sum_{\tilde I\geq\tilde J}
{I\choose J}\cdot{\tilde I\choose\tilde J}a_I\cdot \tilde a_{\tilde I}
(q-p)^{I-J,\tilde I-\tilde J}.$$
Here, the first sum is over all $J,\tilde J$ such that $L=(J,\tilde J)$.
But $\sum_Lc_L(x-q)^L=(\sum_Ja_J(x-q)^J)\cdot(\sum_{\tilde J}
\tilde a_J(x-q)^{\tilde J})$, which is just the product of
$\alpha_r(p,q)(f)$ and $\alpha_r(p,q)(\tilde f)$.
\end{proof}

\begin{beisp}
Let $n=2$, $p=0$ and $f=x_1x_2-tx_2x_1$.
For any $q=(q_1,q_2)\in\CC^2$, 
we get $a_0(q)=q_1q_2-tq_2q_1=(1-t)q_1q_2$.
$a_{(1)}(q)=a_{(1,2)}q^{(2)}+a_{(2,1)}q^{(2)}=(1-t)q_2$,
$a_{(2)}(q)=a_{(1,2)}q^{(1)}+a_{(2,1)}q^{(1)}=(1-t)q_1$,
$a_{(1,2)}(q)=a_{(1,2)}=1$, $a_{(2,1)}(q)=a_{(2,1)}=-t$.
Thus, $\sum_Ja_J(q)(x-q)^J=(1-t)q_1q_2+(1-t)q_2(x_1-q_1)+
(1-t)q_1(x_2-q_2)+(x_1-q_1)(x_2-q_2)-t(x_2-q_2)(x_1-q_1)$.
\end{beisp}

\paragraph{The affine analytic NC space}
For open subsets $U\sub\CC^n$, set $\Oh(U)$ to be the
set of all elements $(f_p)_{p\in U}$ of the product
$\prod_{p\in U}\CC\{x-p\}$ such that the following condition holds:
Let $p$ be a point in $U$ and suppose that $f_p$ converges on
the open polydisk $P(p,r)\sub U$. Then, for each
$q$ in $P(p,r)$, we have $f_q=\alpha_r(p,q)(f_p)$.
It follows by Proposition~\ref{alpha}, that $\Oh(U)$ is
an associative ring. In consequence, the assignment
$$\Oh:U\mapsto\Oh(U)$$ defines
a sheaf of associative $\CC$-algebras on $\CC^n$.
Write $\Oh_{\ab}$ for the sheaf of holomorphic
functions on $\CC^n$. 

\begin{prop}
We have a natural epimorphism
$\ab:\Oh\nach\Oh_{\ab}$ of sheaves of rings, which splits
as a morphism of sheaves of vector spaces.
\end{prop}
\begin{proof}
For open subsets $U$ of $\CC^n$,
we identify $\Oh_{\ab}(U)$ with set of families
$(g_p)_{p\in U}$ of commuative power series
$g_p\in\CC\{x_1-p_1\ddd x_n-p_n\}$, such that the functions
represented by $g_p$ and $g_q$ coincide on the intersection
of the convergence regions of $g_p$ and $g_q$ (if $g_p$ converges
absolutely at $q$ this is equivalent to $\alpha_r(p,q)(g_{p,\unab})=
g_{q,\unab}$). Now, $\ab(U)$ maps an element $(f_p)_{p\in U}$ of
$\Oh(U)$ to the family $(f_{p,\ab})$.
We even get that $\Oh(U)$ is surjective, for any open
subset of $\CC^n$, since $\Oh(U)$ splits by the map
$g=(g_p)_{p\in U}\mapsto(g_{p,\unab})_{p\in U}$. 
\end{proof}


By the following proposition, together with Proposition~\ref{localring},
the pair $(\CC^n,\Oh)$ is a locally associative ringed space.
\begin{prop}\label{lokal}
For each point $p\in\CC^n$, the stalk $\Oh_p$ is isomorphic
to the noncommutative convergent power series ring
$\CC\{x-p\}$.
\end{prop}
\begin{proof}
Exercise.
\end{proof}

\paragraph{Local models}
Now we are able to define local models of NC analytic spaces.
Let $U$ be an open subset of $\CC^n$. Let $\Oh_U$ be
the restriction of $\Oh$ to $U$. 
Let $\I$ be a sheaf of twosided ideals of $\Oh_U$,
such that $\I_{\ab}$ is a coherent $\Oh_{U,\ab}$-module.
Let $V$ be the subset of $U$ given by the zeros
of $\I_{\ab}$.
Write $\iota$ for the inclusion $V\nach U$.
Set $\Oh_V:=\iota^{-1}(\Oh_U/\I)$. 

\begin{prop}
For $q\in V$, the stalk $\Oh_{V,q}$ is naturally isomorphic
to the quotient of $\CC\{x-q\}$ by the twosided
ideal by $\I_q$.
\end{prop}

In consequence,
the pair $V(\I):=(V,\Oh_V)$ is a locally associative ringed space.
Pairs of this form will be called \textbf{local models of
NC complex analytic spaces}. 
Observe that the abelization
$(V,\Oh_{V,\ab})$ is a local model of a complex analytic
space in the classical sense.

\begin{beisp}
Let $\K$ be the kernel of the sheaf homomorphism
$\ab:\Oh_U\nach\Oh_{U,\ab}$. Then $\K_{\ab}=0$ is coherent
and the pair $V(\I)$ is just a local model for a (commutative) analytic
manifold. 
\end{beisp}

\paragraph{NC complex spaces}

Now, we come to the main definition.

\begin{defi}
An \textbf{NC complex analytic space} is a locally associative
ringed space $(X,\Oh_X)$, locally isomorphic to a local model.
An \textbf{NC complex analytic manifold} is a locally associative
ringed space, locally isomorphic to an open subspace
of $\CC^n$ with the canonical NC structure.
\end{defi}

\section{Examples}

\begin{beisp}(Local models)\\
Choose any splitting $s:\Oh_{\ab}\nach\Oh$ of the abelization map
in the category of sheaves of $\CC$-vector spaces.
For any local model $X_{\ab}=V(U,\I_{\ab})$ of a complex analytic space,
set $\I$ to be the two-sided ideal $(s\I)^\wedge$, generated by the
image of $\I_{\ab}$ under $s$. (To generate more examples, add
subsheaves of the commutator sheaf to $\I$.)
We have $\ab(\I)=\I_{\ab}$. Then the NC local model
$X=V(U,\I)$ has the same underlying topological space as
$X_{\ab}$ and $\Oh_{X_{\ab}}=\ab(\Oh_X)$.
\end{beisp}

\begin{beisp}(Analytic supermanifolds)\\
By definition, an analytic supermanifold is a locally ringed space
$(M,\Oh_M)$ with a $\ZZ/2$-graded structure sheaf 
$\Oh_M=\Oh_M^0\oplus\Oh_M^1$ such that the following two conditions hold:
\begin{enumerate}
\item
Let $\Oh_M^{\red}:=\Oh_M/(\Oh_M^1)$ be the quotient of $\Oh_M$ modulo
the two-sided ideal, generated by $\Oh_M^1$. The pair
$(M,\Oh_M^{\red})$ is a usual analytic manifold.
\item
There is an $r\geq 0$ and
each point of $M$ has an open  neighbourhood $P$ in $M$ such
that the restriction $\Oh_M|_P$ of the structure sheaf is
isomorphic to the exterior $\Oh^{\red}_M|_P$-algebra over the
free $\Oh_M^{\red}|_P$-module of rank $r$.
\end{enumerate}
\noindent
Each analytic supermanifold is an NC complex analytic space:
Let $M^{\red}=(M,\Oh_M^{\red})$ be n-dimensional. Choose an
analytic atlas of $M^{\red}$
of the form
$$\{\phi_k:P(p^{(k)},r^{(k)})\nach P^{(k)}|\;k\in K\},$$
where $(P^{(k)})_{k\in K}$ is an open covering of $M$ such
that the restriction $\Oh_M|_{P^{(k)}}$ satisfies
condition (2). Consider the sheaf $\I$ of two-sided ideals
of the NC structure sheaf $\Oh_{P(p^{(k)},r^{(k)})\times\CC^r}$
of $P(p^{(k)},r^{(k)})\times\CC^r$, generated by the following relations:

\[ \begin{array}{rcl}
\lbrack x_i-p^{(k)}_i,x_j-p^{(k)}_j\rbrack & \text{ for }& 1\leq i<j\leq n,\\
\lbrack x_i-p^{(k)}_i,y_j\rbrack & 
\text{ for } & i=1\ddd n\text{ and }j=1\ddd r,\\
y_iy_j+y_jy_i & \text{ for } & 1\leq i\leq j\leq r.
\end{array}\]
Observe that the abelization $\I_{\ab}$ is coherent, and we have
an isomorphism
$$(P^{(k)},\Oh_M|_{P^{(k)}})\isom V(P(p^{(k)},r^{(k)}),\I).$$
\end{beisp}

\begin{beisp}(NC projective space)\\ 
The (commutative) projective space $\CC P^2$ is obtained by glueing
three copies $U_1,U_1$ and $U_2$ of $\CC^2$ along the
following isomorphisms between open subsets:
\[\begin{array}{c}
U_{01}\nach U_{10}, (x_1,x_2)\mapsto(\frac{1}{x_1},\frac{x_2}{x_1}),\\
U_{02}\nach U_{20}, (x_1,x_2)\mapsto(\frac{1}{x_2},\frac{x_1}{x_2}),\\
U_{12}\nach U_{21}, (y_1,y_2)\mapsto(\frac{y_1}{y_2},\frac{1}{y_2}),
\end{array}\]
where $U_{ij}$ is the subset of $U_i$, where the indicated map
is defined.
For each point $p_1\neq 0$ in $\CC$, the fraction $\frac{1}{x_1}$
can be interpreted as a power series in $x_1-p_1$, converging 
near the point $p_1$. In order to use the same identifications
to glue the noncommutative affine spaces $U_0,U_1$ and $U_2$,
we have to choose, if, for $(p_1,p_2)$ in $\CC^2$
with $P_1\neq 0$, we interpret $\frac{x_2}{x_1}$
as the noncommutative power series
$x_2\cdot\frac{1}{x_1}$ in $x_1-p_1,x_2-p_2$ or 
as the noncommutative power series
$\frac{1}{x_1}\cdot x_2$. If, by convention, we always use
the second interpretation, we get compatibility of the gluing maps.
To be precise, the isomorphism 
$\phi_{10}:(U_{01},Oh_{U_{01}}\nach (U_{10},\Oh_{U_{10}}$
of NC complex spaces is given as follows: The map of
topolological spaces $U_{01}\nach U_{10}$ is as in the commutative
case. The associated map $\Oh_{U_{10}}\nach\phi_{10,\ast}\Oh_{U_{01}}$
is stalkwise given by the map $\CC\{y-q\}\nach\CC\{x-p\}$,
sending $y_1-q_1$ to
the power series $\sum_Ia_{1,I}(x-p)^I$, where
$$a_{1,I}=\left\{ \begin{array}{lcl}
\frac{-1}{p_1^{n+1}} & for & I=(1\ddd 1) \quad (n\text{copies})\\
0 & else &
\end{array}\right. $$
and sending $y_2-q_2$ to
the power series $\sum_Ia_{2,I}(x-p)^I$, where
$$a_{2,I}=\left\{ \begin{array}{lcl}
\frac{-1}{p_1^{n+1}}\cdot p_2 & for & I=(1\ddd 1) \quad (n\text{ copies})\\
\frac{-1}{p_1^{n+1}} & for & I=(1\ddd 1,2) \quad (n\text{ copies of }1)\\
0 & else &
\end{array}\right. $$
We leave it to the reader to give an explicit description of the other glueing
isomorphisms.
In this way, we define an NC structure on the topological space 
$\CC P^2$. Exactly in the same way, we get an NC structure
$\Oh_{\CC P^n}$ on $\CC P^n$, for each $n\geq 2$.
The NC complex manifold $(\CC P^n,\Oh_{\CC P^n})$
is called \textbf{NC projective space}. We shall denote it
by $\CC P^n_{NC}$.
We recover the ordinary, commutative projective space as
abelization of the NC complex space. 
\end{beisp}

\begin{beisp}(NC projective varieties)\\
By Chow's theorem, each complex analytic subspace $X$ of
$\CC P^n$ is algebraic, i.e. given by
homogeneous polynomials $f_1\ddd f_m$ in $\CC[x_0\ddd x_n]$.
For each $i=1\ddd m$, we can choose a homogeneous
lift of $f_i$ in $\tilde f_i\in\CC[x_0|...|x_n]$ 
and consider the two-sided ideal $I$ of
$\CC[x_0|...|x_n]$ generated by those lifts.\\
...
\end{beisp}

\end{document}